\documentclass [12pt]{article}
\usepackage{amsmath}
\usepackage{amsfonts, amssymb}
\usepackage[cp1251]{inputenc} 

\newtheorem{thm}{Theorem}
\newtheorem{lemma}[thm]{Lemma}
\newtheorem{prop}[thm]{Proposition}

\newtheorem{cor}[thm]{Corollary}

\def\r{{\mathbb R}}
\def\e{{\mathbb E}}

\def\be#1\ee{\begin{equation}#1\end{equation}}
\newcommand{\bea}{\begin{eqnarray}}
\newcommand{\eea}{\end{eqnarray}}
\newcommand{\beaa}{\begin{eqnarray*}}
\newcommand{\eeaa}{\end{eqnarray*}}

\def\ph{{\varphi}}
\def\X{\mathcal X}

\begin{document}

\title{On Haar Expansion of Riemann--Liouville process
\\
in a critical case
\footnote{The work supported by the RFBR-DFG grant 09-01-91331 "Geometry and asymptotics of random structures".}
}
\author
{M.A. Lifshits}
\date{ }
\maketitle

\begin{abstract}
\noindent We show that Haar-based series representation of the critical
Riemann--Liouville process $R^{\alpha}$ with $\alpha=3/2$
is rearrangement non-optimal in the sense of convergence rate in ${\mathbf C}[0,1]$.
\end{abstract}


\section*{Introduction and main result}
The aim of this note is to solve a problem stated by A.Ayache and W. Linde in their
recent work \cite{AyLin}. In their article, the quality of several common series representations
of fractional processes is considered. In particular, Ayache and Linde examine the representations
of Riemann--Liouville processes $R^{\alpha}$, $\alpha>1/2$, based on Haar and trigonometric systems. 
They were able to show that Haar-based representation of  $R^{\alpha}$ is optimal w.r.t. the uniform norm
when $1/2<\alpha<3/2$ while for $\alpha>3/2$ it is not optimal. Their approach however does not yield
an answer in the delicate critical case $\alpha=3/2$. We will show that Haar-based representation is
not optimal for $R^{3/2}$ either.
\medskip

{\it It is an immense pleasure for the author to stress that the main tool used in this note is due
to V.N.Sudakov whose anniversary we celebrate in this volume.} 
\medskip

The reader can find recent results on series representations of fractional processes and fields in
\cite{AyTaq}, \cite{DZ}, \cite{GilSot}, \cite{Igl}, \cite{KuLi}, \cite{Mal},  \cite{Sch}. 

Recall that Riemann-Liouville process $R^\alpha$ is defined by white noise representation
\[
R^\alpha(t)=\frac{1}{\Gamma(\alpha)} \int_0^t (t-u)^{\alpha-1} dW(u),
\qquad 0\le t\le 1.
\]
This process (RL for short) is known to have continuous sample paths whenever $\alpha>1/2$. 
It belongs to the family of so called fractional processes along with more
widely known  fractional Brownian motion $W^H$ (whenever $\alpha\in (1/2,3/2)$ and $H=\alpha-1/2$) 
and periodic stationary Weil process $I^\alpha$.
All these processes  differ by very smooth terms and therefore their approximation properties 
we discuss here are essentially the same, see \cite{AyLin} for details.
We also refer to  \cite{LifSim} for further properties and stable extensions of RL-process.

The representation of RL-process with $\alpha>1/2$ generated by Haar base writes as follows,
\be \label{haarrepr}
   R^\alpha(t)= \xi_{-1} \frac{t^\alpha}{\Gamma(\alpha+1)}
   +\sum_{j=0}^\infty \sum_{k=0}^{2^j-1} \xi_{j,k}\ (R_\alpha h_{j,k}) (t)
\ee
where $R_\alpha$ is the classical RL integration operator,
\[
R_\alpha h(t)=\frac{1}{\Gamma(\alpha)} \int_0^t (t-u)^{\alpha-1} h(u) du,
\]
$\{\xi_{-1}, (\xi_{j,k}) \}$ is a family of i.i.d. standard normal random variables and $h_{j,k}$'s are
the Haar functions
\[
h_{j,k}(t)= 2^{j/2} \left\{
{\bf 1}_{\left[\frac{2k}{2^{j+1}}, \frac{2k+1}{2^{j+1}}\right)} (t) -
{\bf 1}_{\left[\frac{2k+1}{2^{j+1}}, \frac{2k+2}{2^{j+1}}\right)} (t) 
\right\}.
\]
The series converges almost surely uniformly on $t\in[0,1]$, i.e. it converges in the sense of the uniform
norm $||\cdot||_\infty$.

Recall that one can evaluate the integrated Haar functions by the formula
\be \label{inthaar}
(R_\alpha h_{j,k}) (t) =
\frac{2^{j/2}}{\Gamma(\alpha+1)}
\left\{
\left(t-\frac{2k+2}{2^{j+1}}\right)_+^\alpha -2\left(t-\frac{2k+1}{2^{j+1}}\right)_+^\alpha
+\left(t-\frac{2k}{2^{j+1}}\right)_+^\alpha
\right\}.
\ee

Now we recall the necessary standard notation as well as the 
notions related to finite rank approximation of Gaussian random functions.

Throughout the article  $f_n \succeq g_n$ means $\liminf_n f_n/g_n >0 $ while $f_n \asymp g_n$ means that both
$f_n  \succeq g_n$ and $g_n \succeq f_n$ hold.
We write $\#(B)$ for the number of points in a set $B$, while $|T|$ denotes the length of an interval $T$.  Finally,
$c$ denotes unspecified positive and finite constants which can be different in each occurrence.

Let $X$ be a centered Gaussian random element of a normed space $(\X,||\cdot||)$. The
$\ell$-numbers (stochastic approximation numbers) of $X$ are defined as
\[
\ell_n(X) := \inf_{(\xi_i),(\ph_i)} \left\{ \e\left\|X-\sum_{i=1}^{n-1} \xi_i \ph_i\right\| \right\}.
\]
Here the infimum is taken over all families of random variables $(\xi_i)$ and all
finite deterministic subsets $(\ph_i)\subset \X$. A series representation 
\[
X=\sum_{i=1}^\infty \xi_i \, \ph_i
\]
is called {\it optimal} if
\[
\e\left\|\sum_{i=n}^{\infty} \xi_i\, \ph_i\right\| \asymp \ell_n(X)
\qquad \textrm{as}\ n\to \infty.
\]
It is called {\it rearrangement non-optimal}, if it can not be rendered optimal by any permutation
of $(\ph_i)$. In particular, for non-optimal representation the approximation error
$\e\left\|\sum_{i=n}^{\infty} \xi_i\, \ph_i\right\|$ tends to zero slower than optimal rate $\ell_n(X)$.

In the particular case of RL-process and uniform norm $||\cdot||_\infty$, the optimal approximation rate is well 
known, see \cite{KuLi}. Namely, for any $\alpha>1/2$ it is true that
\[
\ell_n(R^\alpha) \asymp n^{-(\alpha-1/2)}\, \sqrt{\ln n}, 
\qquad \textrm{as}\ n\to \infty.
\]
The optimal rate can be attained, for example, by using modified Daubechies wavelet base, see \cite{Sch}.
The remaining question is to understand which representations provide this rate and which ones perform more
poorly.

For the critical case, $\alpha=3/2$, we have
\be \label{optrate32}
\ell_n(R^{3/2}) \asymp n^{-1}\, \sqrt{\ln n}, 
\qquad \textrm{as}\ n\to \infty.
\ee
Our main result shows that Haar-based representation is rearrangement non-optimal.

\begin{thm} \label{t1} 
Let $\{\xi_i \, \ph_i \}$ be any rearrangement of the Haar-based series representation
$(\ref{haarrepr})$. Then
\[
\e\left\|\sum_{i=n}^{\infty} \xi_i\, \ph_i\right\|_\infty  \succeq  n^{-1}\, \ln n,
\qquad \textrm{as}\ n\to \infty.
\] 
\end{thm}

We observe a gap of order $\sqrt{\ln n}$ w.r.t. the optimal rate (\ref{optrate32}).

\section*{Proofs}

From now on, we consider only the critical case $\alpha=3/2$.

Introduce the function $H(t)=(t-2)_+^{3/2}-2 (t-1)_+^{3/2}+ t_+^{3/2}$ which is the unscaled version
of integrated Haar functions (\ref{inthaar}). 
Namely, we have
\be \label{scale1}
(R_{3/2} h_{j,k}) (t) = \frac
{H\left(2^{j+1}\left(t-\frac{2k}{2^{j+1}}\right)\right)}
{2^{3/2+j}\Gamma(5/2)}\ .
\ee

Notice that for $t\ge 2$
\be \label{asympH}
H(t)= \frac 32 \int_{t-1}^t (x^{1/2}-(x-1)^{1/2})dx\sim \frac 34 \ t^{-1/2}, \qquad 
\textrm{as}\quad t\to \infty.
\ee
For any level number $j$ and a set of positive integers $K\subset [0,\dots,2^j-1]$ we consider a Gaussian process 
$\{X_K(t),t\in[0,2^{j+1}]\}$ defined by
\[
X_K(t)= \sum_{k\in K} \xi_k H(t-2k)
= \sum_{k\in K, k\le t/2} \xi_k H(t-2k),
\]
where $\{\xi_k\}$ is a family of i.i.d. standard normal random variables. Clearly, $X_K$ is the unscaled version
of the $K$-part of level $j$ in Haar-based representation (\ref{haarrepr}), since by (\ref{scale1})
\be \label{scale2} 
\sum_{k\in K} \xi_{j,k}\ (R_{3/2} h_{j,k}) (t)
=
\frac{X_K(2^{j+1}t)}
{2^{3/2+j}\Gamma(5/2)} 
\ee
in distribution.
\medskip

\begin{prop} \label{p1} There exists a numerical constant $c$ such that
\[
\e \sup_{t\in[0,2^{j+1}]}|X_K(t)| \ge c\ j,
\]
for all $j>0$ and all $K\subset \{0,\dots,2^j-1\}$ such that $\#(K)\ge 2^{j-1}$.
\end{prop}

\begin{cor}\ By scaling $(\ref{scale2})$ it follows that
\be \label{c1}
\e \sup_{t\in[0,1]}\left|\sum_{k\in K} \xi_{j,k}(R_{3/2} h_{j,k})(t)\right | \ge c\ j \, 2^{-j},
\ee
where $h_{j,k}$ are Haar functions.
\end{cor}

{\bf Proof of Proposition \ref{p1}.}\  
First of all, let us evaluate the natural distance associated to the process $X_K$. 
Let $0\le s<t\le 2^{j+1}$. We have
\begin{eqnarray*}
X_K(t) &=& \sum_{k\in K: k\le t/2} \xi_k H(t-2k) 
\\
 &=& \sum_{k\in K:k\le s/2} \xi_k H(t-2k)+ \sum_{k\in K: s/2< k\le t/2} \xi_k H(t-2k) 
\\
 &:=& Y(s,t)+Z(s,t)\ . 
 \end{eqnarray*}
Notice that the variable $X_K(s)$ belongs to
$span(\xi_k, k\in K, k \le s/2)$ and the same is true for $Y(s,t)$, while $Z(s,t)$ is orthogonal to that span. 
By using (\ref{asympH}), it follows that for $s<t$ we have
\begin{eqnarray*}
\e (X_K(t)-X_K(s))^2 &\ge& \e \, Z(s,t)^2=  \sum_{k\in K: s/2< k\le t/2} H(t-2k)^2  
\\
 &\ge&    c\ \sum_{k\in K: s/2< k\le t/2-1} (t-2k)^{-1}    \, .
 \end{eqnarray*}
 Suppose now that we found in $[0,2^{j+1}]$ an ordered subsequence $t_1,\dots, t_m$,
 $m\ge 2^{j/2-2}$ such that for any $2\le i\le m$ it is true that
 \be \label{wewant}
    \sum_{k\in K: t_{i-1}/2< k\le t_i/2-1} (t_i-2k)^{-1} \ge c\, j.
 \ee
 
Then, by Sudakov lower bound (see e.g. \cite{GRF}, Section 14)
\begin{eqnarray*}
\e \sup_{t\in[0,2^{j+1}]}|X_K(t)| &\ge& \e \sup_{1\le i \le m} X_K(t_i) 
\\
&\ge& c (\ln m)^{1/2} \cdot \inf_{i_1\not =i_2} 
\left[  \e(X(t_{i_1})-X(t_{i_2}))^2 \right]^{1/2}
\\ 
&\ge&  c\, j^{1/2}\cdot c\, j^{1/2} = c\ j\, .
\qquad
\end{eqnarray*}  
\bigskip
It remains to find a sequence $(t_i)$ satisfying assumption (\ref{wewant}). 
Would $K$ coincide with the maximal set, $K=\{0,\dots 2^j-1\}$,  we could just take
$t_i= i 2^{j/2}, i=0,\dots, 2^{j/2}-1$. However, the situation is more delicate for general case. 
We will use the following elementary fact worth to be stated separately.

\begin{lemma} \label{l1} Let the points
$s_1 < s_2 < \dots < s_{2q}$ belong to an interval $T$ in $\r$. Then
\[
\max_{1 < r\le 2q} \sum_{i=1}^{r-1} \frac{1}{s_r-s_i} \ge c\ \frac q{|T|} \ \ln q.
\]
\end{lemma}

{\bf Proof of Lemma \ref{l1}}. Take an integer $j$ such that $1\le j\le q$.
We have
\begin{eqnarray*}
 \sum_{r=q+1}^{2q} (s_r-s_{r-j}) &=& \sum_{r=q+1}^{2q} \sum_{l=r-j}^{r-1} (s_{l+1}-s_{l})
\\
&\le&   \sum_{l=q+1-j}^{2q-1} j (s_{l+1}-s_{l}) \le j |T|.
\end{eqnarray*}
By convexity of the function $x\to \frac 1x$ we have
\[
 \frac 1q \sum_{r=q+1}^{2q}
  \frac 1{s_r-s_{r-j}}   
  \ge \left( \frac 1q  \sum_{r=q+1}^{2q} (s_r-s_{r-j})\right)^{-1} 
  \ge \left( \frac{j|T|}{q} \right)^{-1}.
\]
By summing up over $j$ we get
\begin{eqnarray*}
\frac{1}{q} \sum_{r=q+1}^{2q}\sum_{j=1}^{q} \frac 1{s_r-s_{r-j}}
&=&
 \sum_{j=1}^{q} \frac{1}{q}  \sum_{r=q+1}^{2q}  \frac 1{s_r-s_{r-j}}
\\
&\ge&
\sum_{j=1}^{q}  \left( \frac{j|T|}{q} \right)^{-1} =\frac q{|T|} \sum_{j=1}^{q}j^{-1}
\ge \frac{cq}{|T|} \ \ln q.
\end{eqnarray*}

Clearly,
\begin{eqnarray*}
\max_{1< r\le 2q} \sum_{i=1}^{r-1} \frac{1}{s_r-s_i} 
&\ge&
\max_{q+1\le r\le 2q} \sum_{i=r-q}^{r-1} \frac{1}{s_r-s_i}
= \max_{q+1\le r\le 2q} \sum_{j=1}^{q} \frac{1}{s_r-s_{r-j}}
\\
&\ge&
\frac{1}{q} \sum_{r=q+1}^{2q}\sum_{j=1}^{q} \frac 1{s_r-s_{r-j}}
\ge \frac{cq}{|T|} \ \ln q. \qquad \Box
\end{eqnarray*}

We continue the proof of Proposition \ref{p1}.
For the sake of expression simplicity, assume that $j$ is an even number, thus $2^{j/2}$
is an integer. We split the interval $[0,2^{j}-1]$ in $2^{j/2}$ blocks
$B_\iota=\left[\iota2^{j/2},(\iota+1)2^{j/2}\right)$.

Let
\[
I=\left\{ \iota:\#(K\cap B_\iota)\ge \frac 14\ 2^{j/2}
\right\}
\]
and $m=\#(I)$. We show that $m$ is large enough, since
\begin{eqnarray*}
2^{j-1}&\le& \#(K) = \sum_{\iota=0}^{2^{j/2}-1}\#(K\cap B_\iota)
\\
&=& \sum_{\iota\in I} \#(K\cap B_\iota) + \sum_{\iota\not\in I}\#(K\cap B_\iota)
\\
&\le& m 2^{j/2} + \left(2^{j/2}-m\right) \frac 14\ 2^{j/2}
\end{eqnarray*}
yields
\[
\frac 12\le \frac {m}{2^{j/2}} +\frac 14 -  \frac 14 \ \frac {m}{2^{j/2}},
\]
thus
\[
m\ge \frac 13\ 2^{j/2}> 2^{j/2-2}.
\]
For each block $B_\iota, \iota\in I$, we apply Lemma \ref{l1} with $T=B_\iota$, $|T|=2^{j/2}$ and
$\{s_1, \dots, s_{2q}\} \subset K\bigcap B_\iota$ and choosing $q$ as large as possible,
\[
2q\ge \#\left(K\cap B_\iota\right)-1 \ge \frac 14\ 2^{j/2} -1.
\]
By Lemma \ref{l1}, we can find an integer point $v_\iota\in \{s_1, \dots, s_{2q}\} \subset K\bigcap B_\iota$ such that
\[
\sum_{s\in K\bigcap B_\iota, s< v_\iota}  \frac{1}{v_\iota-s}  \ge \frac{cq}{T} \ \ln q \ge c\, j.
\]
By renumbering the points $\{2v_\iota,\iota\in I\}$ we get a large collection of
$m\ge 2^{j/2-2}$ points $t_i$ with the property
\begin{eqnarray*}
 \sum_{k\in K: t_{i-1}/2< k\le t_i/2-1} (t_i-2k)^{-1} &=&
 \frac 12 \sum_{k\in K: t_{i-1}/2< k\le t_i/2-1} (t_i/2-k)^{-1}
 \\
 &\ge& \frac 12 \sum_{k\in K: t_{i-1}/2< k\le t_i/2} (t_i/2-k)^{-1} -1
  \\
 &\ge& \frac 12 \sum_{k\in K\bigcap B_\iota, k< v_\iota} (v_\iota-k)^{-1} -1
 \ge cj,
\end{eqnarray*}
where the index $\iota$ is such that $t_i/2=v_\iota$. We also used here the fact that all $v_\iota$'s and all $k$'s
are integers, thus we don't loose more than one while dropping one $k$ from the sum.
Therefore, the collection $\{t_i, i\le m\}$ satisfies (\ref{wewant}) and 
Proposition \ref{p1} is proved. $\Box$
\medskip 

{\bf Proof of Theorem \ref{t1}.} Actually our theorem follows from Proposition \ref{p1} and Lemma 2.3
in \cite{AyLin} immediately. However, we recall the proof for reader's convenience.

Let $\{\xi_i \, \ph_i \}$ be any rearrangement of the Haar-based series representation
$(\ref{haarrepr})$. Take an integer $j\ge 1$, let $n=2^{j-1}$ and let $K\subset [0,\dots,2^j-1]$
be the set of all $k$ such that $R_\alpha h_{j,k}$ belongs to the set $\{\ph_i, i\ge n\}$.
Clearly, $\#(K)\ge 2^{j-1}$, and we can write
\[
\sum_{i\ge n} \xi_i \ph_i = \sum_{k\in K}\xi_{j,k}\ (R_\alpha h_{j,k}) + Y,
\]
where $Y$ and the sum over $K$ are independent. By standard arguments based on Anderson inequality 
for Gaussian processes (see \cite{GRF}), we obtain
\[
\e \left\|\sum_{i\ge n} \xi_i \ph_i\right\|_\infty \ge \e \left\|\sum_{k\in K}\xi_{j,k}\ (R_\alpha h_{j,k})
\right\|_\infty .
\]
Now (\ref{c1}) yields
\[
\e \left\|\sum_{i\ge n} \xi_i \ph_i\right\|_\infty \ge
c\ j \, 2^{-j} \ \asymp n^{-1} \ln n,
\]
and we are done. $\Box$
\medskip 

The author is grateful to A.Ayache and to W.Linde for interesting discussions and for 
attracting his attention to the problem solved in this note.  
\bigskip

\bibliographystyle{amsplain}

\end{document}